\documentclass{amsart}
\usepackage{amssymb,amsmath,amsthm,mathtools, geometry, amssymb, verbatim, esint, soul}
\usepackage{a4wide}
\usepackage{float}
\usepackage{graphicx}
\usepackage[utf8]{inputenc} 
\usepackage{a4wide}
\usepackage{tikz}
\usetikzlibrary{calc}
\usepackage{ulem}
\usepackage{stackengine,scalerel,graphicx}
\stackMath
\usepackage{hyperref}
\usepackage{amsfonts} 
\usepackage{latexsym}
\usepackage[font=small,format=hang,labelfont={sf,bf}]{caption}
\usepackage{epsfig}
\usepackage{subfig}
\usepackage{url}
\usepackage{varioref}
\usepackage{bm}
\usepackage{color}
\usepackage{mathrsfs}
\usepackage{latexsym}
\usepackage{bbm}
\usepackage{faktor}
\usepackage{yhmath}
\usepackage{empheq}
\usepackage{mathtools}
\usepackage{color}
\usepackage{marginnote}
\usepackage{comment}
\usepackage{cancel}
\usetikzlibrary{arrows.meta}

  \newtheorem{theorem}{Theorem}[section]
\newtheorem*{theorem*}{Theorem}
\newtheorem{lemma}[theorem]{Lemma}
\newtheorem{proposition}[theorem]{Proposition}
\newtheorem{corollary}[theorem]{Corollary}

\theoremstyle{definition}
\newtheorem{definition}[theorem]{Definition}
\newtheorem{remark}[theorem]{Remark}

\makeatletter

\@addtoreset{equation}{section}
\makeatother

\newcommand{\EEE}{\color{black}}
\newcommand{\R}{\mathbb{R}}
\newcommand{\Rd}{\mathbb{R}^d}
\newcommand{\N}{\mathbb{N}}
\newcommand{\Z}{\mathbb{Z}}
\newcommand{\dx}{\, \mathrm{d}x}

\newcommand{\norm}[1]{\left\lVert #1\right\rVert}

\newcommand{\Permd}{{\mathcal P}_{d}}
\newcommand{\card}{\#}
\date{}

\title{The quantitative Faber-Krahn inequality for the combinatorial Laplacian in $\mathbb{Z}^d$}

\author[M.~Cicalese]{Marco Cicalese}
\address[Marco Cicalese]{Technical University of Munich, School of Computation Information and Technology, Department of Mathematics, Boltzmannstraße 3,
85748 Garching, Germany}
\email{cicalese@ma.tum.de}

\author[L.~Kreutz]{Leonard Kreutz}
\address[Leonard Kreutz]{Technical University of Munich, School of Computation Information and Technology, Department of Mathematics, Boltzmannstraße 3,
85748 Garching, Germany}
\email{kleo@cit.tum.de}

\author[G.P.~Leonardi]{Gian Paolo Leonardi}
\address[Gian Paolo Leonardi]{Department of Mathematics, University of Trento, 38050 Povo (TN), Italy}
\email{gianpaolo.leonardi@unitn.it}

\author[G.~Morselli]{Gabriele Morselli}
\address[Gabriele Morselli]{Department of Mathematics, University of Trento, 38050 Povo (TN), Italy}
\email{gabriele.morselli@unitn.it}

\begin{document}

\begin{abstract}
While the classical Faber-Krahn inequality shows that the ball uniquely minimizes the first Dirichlet eigenvalue of the Laplacian in the continuum, this rigidity may fail in the discrete setting. We establish quantitative fluctuation estimates for the first Dirichlet eigenvalue of the combinatorial Laplacian on subsets of \(\mathbb{Z}^d\) when their cardinality diverges. Our approach is based on a controlled discrete-to-continuum extension of the associated variational problem and the quantitative Faber-Krahn inequality.
\end{abstract}

\maketitle

\author{}
  

\section{Introduction}

Many classical problems in spectral theory and geometric analysis focus on identifying extremal domains for various functionals. A prominent example is the \textit{Faber–Krahn inequality}, which asserts that among all open sets \( \Omega \subset \mathbb{R}^d \) of fixed volume, Euclidean balls uniquely minimize the first Dirichlet eigenvalue of the Laplacian. In precise terms, if
\[
\lambda(\Omega) = \min\Biggl\{ \int_{\Omega} |\nabla u(x)|^2\,\mathrm{d}x \colon u\in H^1_0(\Omega),\,\|u\|_{L^2(\Omega)}=1 \Biggr\}\,,
\]
then one has
\begin{equation}\label{ineq:Faber_Krahn}
\lambda(\Omega) \ge \lambda(B_{r_{|\Omega|}})\,,
\end{equation}
where $B_{r_{|\Omega|}}\subset \mathbb{R}^d$ is a ball with the same measure as $\Omega$. Moreover the equality holds if and only if \(\Omega\) agrees with $B_{r_{|\Omega|}}$ up to translations and sets of zero $2$-capacity. Even though in Euclidean space minimizers of the first eigenvalue are rigid, rigidity may fail in more general settings. This is also the case of graphs.

Spectral graph theory (see for instance \cite{BLS,Grigoryan}) has attracted the interest of a wide mathematical community as it reveals remarkable connections with differential geometry, Riemannian geometry, algebraic graph theory, and probability theory (see \cite{Biggs,Chung} and references therein). Furthermore, it has applications in other areas of science such as theoretical chemistry \cite{Tompa}, computer science, and physics \cite{Kurasov}. While for regular trees it has been proved in \cite{Pruss} that minimizers of the first eigenvalue of the combinatorial Laplacian are essentially rigid (up to graph automorphism), this is not true in general for Bravais lattices. In this paper, we will consider the special case of the integer lattice $\mathbb{Z}^d$. Here, one defines the discrete first Dirichlet eigenvalue of the combinatorial Laplacian for a finite set $ X\subset \mathbb{Z}^d$ (with cardinality $\#X=N$) by
\begin{equation}\label{intro:lambda}
\lambda_N\big(X\big) = \min\Biggl\{N^{\frac{2-d}{d}} \sum_{\substack{i,j\in\mathbb{Z}^d \\ |i-j|=1}} \big|u(i)-u(j)\big|^2 \colon u(i)=0 \text{ in } \mathbb{Z}^d\setminus X,\; \frac{1}{N}\sum_{i\in\mathbb{Z}^d} u(i)^2 = 1 \Biggr\}\,.
\end{equation}
By the discrete to continuum $\Gamma$-convergence result in \cite{Ali-Cic}, as $N\to \infty$ the discrete first Dirichlet eigenvalue of the combinatorial Laplacian on $\mathbb{Z}^d$ converges to the continuum one. Consequently, one might expect that, under a volume constraint, minimizing sets of $\lambda_N$ converge to the Euclidean ball (which minimizes the continuum problem), as proven in \cite{S-R}. However, there is no precise characterization of minimal sets for fixed $N$ and one might even expect that minimizers of the discrete problem are generically not unique. More precisely, there might be arbitrarily large $N$ for which there exist distinct subsets \( X_N, Y_N\subset \mathbb{Z}^d \) such that
\begin{equation}\label{intro:min}
\lambda_N(X_N) = \lambda_N(Y_N) = \min\Bigl\{ \lambda_N(X) \colon X\subset \mathbb{Z}^d,\; \#X=N \Bigr\}\,,
\end{equation}
yet for every discrete isometry \(R\colon\mathbb{Z}^d\to \mathbb{Z}^d\), $R(x)=Qx+\tau$, with $Q \in O(d,\mathbb{Z})$ and $\tau \in \mathbb{Z}^d$, one has
\[
\#\Bigl( X_N\Delta R(Y_N) \Bigr) \neq 0\,.
\]
Note that this property is in contrast to the characterization of minimizers for the regular trees as the discrete isometries are the automorphism group of $\mathbb{Z}^d$. As it is customary for geometric problems on lattices, the lack of uniqueness calls for \textit{fluctuation estimates}, that is, estimates quantifying the deviation of different minimizers from each other as \( N \) grows. Even more, one can obtain such fluctuation estimates for almost minimizers of the problem under consideration.

In this paper, we establish, for the first time in the spectral setting, discrete fluctuation estimates for the first Dirichlet eigenvalue of the combinatorial Laplacian on $\mathbb{Z}^d$. In its simplest formulation, our main result (see Theorem~\ref{main_THM}) shows that if \( X_N, Y_N\subset \mathbb{Z}^d \) satisfy \eqref{intro:min}, then there exists a constant \( C_d>0 \), depending only on the dimension, such that
\begin{equation}\label{intro:fluc}
\#\Bigl( X_N\Delta Y_N \Bigr) \le C_d\, N^{1-\frac{1}{2d}}\,.
\end{equation}
We note that in Theorem~\ref{main_THM} we prove more general estimates than the one above. They hold for almost minimizers provided the energy gap from the minimum is not too big (and some non-uniform control on the perimeter of the sets is assumed). Fluctuation inequalities like the one in \eqref{intro:fluc} have been the object of recent studies (see \cite{CL,Davoli-Piovano-Stefanelli16,Davoli-Piovano-Stefanelli17,FriedrichKreutz:19,FriedrichKreutz:20, Mainini-Piovano-Schmidt-Stefanelli2019,Mainini-Piovano-Stefanelli2019,Schmidt:13}). In particular, in \cite{CL} the authors proved optimal fluctuations estimates for the edge isoperimetric problem (see \cite{Bollobas,CKL}) in the planar setting by establishing a connection between fluctuation estimates and the quantitative isoperimetric inequalty in the continuum setting. In Section~\ref{sec:fluctuations}  of this paper we follow the same strategy, this time substituting the quantitative isoperimetric inequality \cite{Fusco-Maggi-Pratelli} with the the quantitative Faber-Krahn inequality from \cite{BDPV}. The latter states that there exists a constant $C_d>0$, depending only on the dimension, such that for any measurable set $\Omega\subset\mathbb{R}^d$ with finite measure,
\[
\inf_{z\in\mathbb{R}^d} \big|\Omega\Delta (B_{r_{|\Omega|}}+z)\big| \le C_d\, |\Omega|\left(\frac{\lambda(\Omega)-\lambda(B_{r_{|\Omega|}})}{\lambda(B_{r_{|\Omega|}})}\right)^{1/2}\,.
\]
We outline our strategy shortly below. Roughly speaking, a fluctuation estimate like the one in \eqref{intro:fluc} can be obtained from the inequality above by a  two steps procedure. In a first step we associate a continuum measurable set $\zeta(X)$ to a discrete configuration $X$ in such a way that $|X\Delta(z+B_{r_{|\Omega|}}\cap\mathbb{Z}^d)|\simeq |\zeta(X)\Delta(z+B_{r_{|\Omega|}})|$. In a second step we derive optimal upper bounds on $\lambda(\zeta(X))-\lambda_N(X)$. As the first eigenvalue of the Dirichlet Laplacian on a given set is itself obtained by minimizing the Reyleigh quotient, the approximation procedure requires also the extension of the function $u$ in the definition of $\lambda_N(X)$ in \eqref{intro:lambda}. Both the extension of the set $X$ and of the function $u$ are obtained making use of the Kuhn decomposition of the cube (see Section~\ref{subsec:Kuhn}).
It is worth mentioning that the error introduced in the extension procedure involves the perimeter of a discrete set with finite first Dirichlet eigenvalue of the Laplacian. One can get rid of that dependence in dimension $d=2$ in the case of minimizers and for $d>2$ in the case of symmetric minimizers, whose existence is guaranteed by a discrete Riesz rearrangement inequality (see Theorem~\ref{thm:discrete_Riesz_rearrangement} and Corollary~\ref{cor:sym}). \\

The remainder of the paper is organized as follows. In Section~\ref{sec:prelim} we introduce the necessary preliminaries and define the discrete first eigenvalue of the combinatorial Laplacian. In Section~\ref{sec:prop_optimal_sets} we discuss elementary properties of optimal sets for the first eigenvalue of the combinatorial Laplacian.  Finally, Section~\ref{sec:fluctuations} is devoted to the derivation of all discrete-to-continuum estimates to establish the main quantitative fluctuation result, thereby highlighting the differences and connections between the discrete and continuum spectral problems.  \EEE

\section{Notation and preliminaries} \label{sec:prelim}

In the following, given $\Omega\subset\Rd$, we denote by $\mathring{\Omega}$ its interior and by $(\Omega)_{r}$ its $r$-neighborhood, namely 
\begin{align*}
    (\Omega)_{r} = \{ x\in \R^d \colon  \mathrm{dist}(x,\Omega) < r\}\,.
\end{align*}  
With a little abuse of notation, the symbol $|\cdot|$ denotes both the Lebesgue measure of a set in $\Rd$ and the standard Euclidean norm of a vector in $\mathbb{R}^d$. For $k=1,\dots,d$ we denote by $e_k$ the $k^{th}$ element of the canonical basis of $\Rd$. For $x\in \mathbb{R}^d$ and $r>0$, $Q_r(x)$ denotes the closed coordinate cube centered at $x$ and with side lengths $r$, while  $B_r(x)$ the closed Euclidean ball of radius $r$ centered at $x$. When no confusion is possible we also use the notation $B_r=B_r(0)$.  We denote by $\chi_{\Omega}$ the characteristic function of $\Omega \subset \mathbb{R}^d$. Given two sets $E,F \subset \mathbb{R}^d$ we denote by $E\Delta F$ their symmetric difference. Given two points $x,y \in \mathbb{R}^d$, we write $[x,y]$ for the closed segment joining $x$ and $y$. We denote by $\mathbb{Z}^d$ the $d$-dimensional integer lattice and by $\mathcal{N}$ the set of pairs of neighboring points of $\mathbb{Z}^d$, namely
\begin{equation*}
    \mathcal{N}=\left\{(p,q)\in\Z^d\times\Z^d \colon |p-q| = 1 \right\}\,.
\end{equation*}
We set $\mathcal{X}$ to be the collection of all subsets $X\subset \mathbb{Z}^d$. For $N \in \mathbb{N}$ we further set $\mathcal{X}_N:= \{ X \in \mathcal{X} \colon \# X=N\}$. 
We denote by $\mathcal{P}_d$ the set of permutations of the set $\{1,\ldots,d\}$. Furthermore, throughout all the following estimates, we will call $C_d$ a constant depending only on the dimension $d$, whose value may change from line to line.

\subsection{Kuhn decomposition}  \label{subsec:Kuhn}
We consider the decomposition of the closed, unit cube $[0,1]^d \subset\Rd$ into $d!$ simplices $\{T_\pi\}_{\pi \in \mathcal{P}_d}$, where
$$T_\pi = \{x\in \R^d\colon 0\le x_{\pi(d)}\le x_{\pi(d-1)}\le\dots\le x_{\pi(1)}\le 1\}\,.$$
We note that $\mathring{T_\pi} \cap  \mathring{T_{\pi'}}=\emptyset$ for $\pi,\pi'\in \mathcal{P}_d$ and $\pi \neq \pi'$. Given $z\in \Z^d$ we set 
\[
T_{\pi}(z) = z+T_{\pi}\,.
\]
We denote by $\mathcal{T}$ the \textit{Kuhn decomposition} of $\R^d$, namely
\begin{align*}
\mathcal{T}=\{T\colon T= T_{\pi}(z) \text{ for some } z \in \mathbb{Z}^d \text{ and } \pi \in \mathcal{P}_d\}.
\end{align*}
Given $z \in \mathbb{Z}^d$, we introduce also the notation
\begin{align*}
\mathcal{T}_z = \{T \in \mathcal{T} \colon T\cap \{z\} \neq \emptyset\}\,.
\end{align*}

\begin{proposition}\label{Kuhn_properties_PROP}
The Kuhn decomposition satisfies the following properties:
    \begin{enumerate}
        \item[{\rm (a)}] $\bigcup_{z \in \mathbb{Z}^d}\mathcal{T}_z = \Rd$  
        \item[{\rm (b)}] For all $z\in \mathbb{Z}^d$ and $\pi \in \mathcal{P}_d$ we have $|T_{\pi}(z)| = 1/d!$;
        \item[{\rm (c)}] For any $k\in \{1,\ldots,d\}$, there exists a unique edge of the simplex $T_{\pi}(z)$ parallel to $e_k$;
        \item[{\rm (d)}] For each $(i,j)\in \mathcal{N}$, there exist $d!$ distinct simplices sharing the segment $[i,j]$ as a common edge;
        \item[{\rm (e)}] For each point $i \in \mathbb{Z}^d$, there exist $(d+1)!$ distinct elements of $\mathcal{T}$ sharing $i$ as a common vertex.
    \end{enumerate}
\end{proposition}

\begin{proof}
Properties \rm{(a)}, \rm{(b)} and \rm{(c)} are a trivial consequence of the definition of $T_{\pi}(z)$. Let us prove \rm{(d)}. 
We preliminarily observe that, given an edge $[i,i+e_{k}]$ with $i\in \Z^{d}$ and $k\in \{1,\dots,d\}$, and a permutation $\pi\in \Permd$, there exists a unique $z\in \Z^{d}$ such that $[i,i+e_{k}]\subset T_{\pi}(z)$. In other words,
\begin{equation}\label{eq:edgeinTzpi}
\card \{z\in \Z^{d}\colon [i,i+e_{k}]\subset T_{\pi}(z)\} = 1\,.
\end{equation} 
Indeed, if $z,z'\in \Z^{d}$ are distinct, then necessarily $T_{\pi}(z)\cap T_{\pi}(z')$ is either empty or contains a single point of the lattice. Then by \eqref{eq:edgeinTzpi} we conclude the proof of \rm{(d)} as
\begin{align*}
\card \{(\pi,z)\in \Permd\times \Z^{d}\colon [i,i+e_{k}]\subset T_{\pi}(z)\} 
= \sum_{\pi\in \Permd}\card \{z\in \Z^{d}\colon [i,i+e_{k}]\subset T_{\pi}(z)\}
= \sum_{\pi\in \Permd} 1 = d!\,.
\end{align*}
To prove \rm{(e)} we first observe that, since $[0,1]^d$ is the periodicity cell of the decomposition, counting the number of different simplices passing through a point is equivalent to summing the number of times each vertex of $[0,1]^d$ belongs to a different simplex of the Kuhn decomposition of $[0,1]^d$. If we call this number $S(d)$, then we have
\begin{align*}
    S(d)& = \sum_{i\in \{0,1\}^d} \card\{\pi\in\mathcal{P}_d \colon i\in T_{\pi}\}= \sum_{i\in \{0,1\}^d} \sum_{\pi\in\Permd} \chi_{T_{\pi}}(i)\\
    & = \sum_{\pi\in\mathcal{P}_d} \sum_{i\in T_{\pi} \cap \{0,1\}^d} 1 = \sum_{\pi\in\mathcal{P}_d} (d+1)=(d+1)!\,.
\end{align*}
\end{proof}

\subsection{Discrete Dirichlet and Perimeter functionals}
In this section, given a discrete set $X\subset\mathbb{Z}^d$ and a function $u:X\to\mathbb{R}$, we define several discrete functionals associated with $X$ and $u$, that will be considered in the rest of the paper. Along with them, we also define their scaled version, corresponding to those energy functionals per unit particle, thus highlighting the dependence on the cardinality of $X$. To this end, given $X\subset \Z^d$ we first define the valence of a point $p\in X$ as
\begin{equation*}
  \mathrm{val}(p)=\card\big\{q\in \mathbb{Z}^d\setminus X \colon (p,q) \in \mathcal{N}\big\}\,.
\end{equation*}
With this definition at hand, we can now introduce the so-called 'edge perimeter' (later on simply perimeter) of a discrete set, as well as its scaled version, as follows.
\begin{definition}\label{def:Per} For $X\in \mathcal{X}$, the discrete perimeter of $X$ is defined as
\[
P(X)=\sum_{p\in X} \mathrm{val}(p)\,.
\]
For $X \in \mathcal{X}_N$ the scaled discrete perimeter of $X$ is defined as
\begin{align*}
P_N(X)=N^{-\frac{d-1}{d}} P(X)\,.
\end{align*}
\end{definition}
\noindent The Dirichlet energy of a scalar function defined on $\mathbb{Z}^d$ is given here below, followed by its scaled version.
\begin{definition} \label{def:EN} Given $u \colon \mathbb{Z}^d \to \mathbb{R}$ we define the discrete Dirichlet energy of $u$ as
\begin{align*}
D(u) = \sum_{(i,j) \in \mathcal{N}} |u(i)-u(j)|^2\,.
\end{align*}  
Given $X\in \mathcal{X}_N$ and $u \colon \mathbb{Z}^d \to \mathbb{R}$ such that $\mathrm{supp}(u)\subseteq X$, we define the scaled discrete Dirichlet energy of $u$ in $X$ as 
\begin{align*}
E_N(u) = N^{-\frac{d-2}{d}} D(u)\,.
\end{align*}
\end{definition}

\begin{remark}[Energy scalings] \label{rem:energyscalings} {\rm The scalings of $E_N$ and $P_N$ are justified by the following $\Gamma$-convergence results. 
\begin{itemize}
    \item[{\rm (a)}] Given a function $u \colon \mathbb{Z}^d \to \mathbb{R}$ we define its piecewise-constant interpolation (subordinated to $N^{-1/d} \mathbb{Z}^d$) as
    \begin{align*}
        \overline{u}_N(x) = u(N^{1/d}z) \quad \text{for } x\in Q_{N^{-1/d}}(z)\,, \quad z \in N^{-1/d} \mathbb{Z}^d\,.
    \end{align*} 
    For $u \colon \mathbb{Z}^d \to \mathbb{R}$ such that
    \begin{align*}
        \frac{1}{N} \sum_{i \in \mathbb{Z}^d} u(i)^p=1
    \end{align*}
    we have $\| \overline{u}_{N}\|_{L^p(\R^d)}=1$. Furthermore,  if $u = 0$ on $\mathbb{Z}^d \setminus N^{\frac{1}{d}} \Omega$ for some $\Omega \subseteq \mathbb{R}^d$ open and bounded, we have that $\mathrm{supp}(\overline{u}_N) \subseteq {(\Omega)_{N^{-1/d}\sqrt{d}}}$ for all $N \in \mathbb{N}$.
    \item[{\rm (b)}] Thanks to the interpolation in {\rm (a)}, with slight abuse of notation one can read $E_N$ as defined on $L^2(\mathbb{R}^d)$. Furthermore, given $\Omega \subset \mathbb{R}^d$ open and bounded, one can introduce the functional
    \begin{align*}
        D_N(v,\Omega) = \begin{cases} 
            E_N(u) &\text{if } v=\overline{u}_N \text{ for } u \colon \mathbb{Z}^d \to \mathbb{R}\,, \frac{1}{N} \sum_{i \in \mathbb{Z}^d} u^2(i)=1\,, u = 0 \text{ on } \mathbb{Z}^d \setminus N^{\frac{1}{d}} \Omega\,, \\
            +\infty &\text{otherwise}\,.
        \end{cases}
    \end{align*}
    For such a functional, thanks to \cite[Theorem~3.1 and Remark~3.2]{Ali-Cic} and the properties of the interpolation described in {\rm (a)}, the following convergence result holds true:
    \begin{align*}
    \Gamma(L^2(\mathbb{R}^d))\text{ -}\lim_{N\to +\infty} D_N(v,\Omega) = \begin{cases} 
        \int_{\mathbb{R}^d}|\nabla v|^2\,\mathrm{d}x &\text{if } v \in H^1(\mathbb{R}^d)\,, v=0 \text{ on } \mathbb{R}^d\setminus \Omega\,, \\& \|v\|_{L^2(\mathbb{R}^d)}=1\,,\\
        +\infty &\text{otherwise on }L^2(\mathbb{R}^d)\,.
        \end{cases}
    \end{align*}
    \item[{\rm (c)}] Arguing as above, again with a slight abuse of notation, one can read $P_N$ as defined on $L^1(\mathbb{R}^d)$ as follows :
    \begin{align*}
        P_N(v) = \begin{cases} P_N(X) &\text{if } v=\overline{u}_N \text{ for } u = \chi_{X} \colon \mathbb{Z}^d \to \{0,1\}\,, \card X = N\,,\\
        +\infty &\text{otherwise.}
        \end{cases}
    \end{align*}
    For such a functional, thanks to \cite[Theorem~4]{Ali-Bra-Cic} and the properties of the interpolation described in {\rm (a)}, the following convergence result holds true:
    \begin{align*}
        \Gamma(L^1(\mathbb{R}^d))\text{ -}\lim_{N\to +\infty} P_N(v) = \begin{cases} \int_{\partial^*E}\|\nu_E\|_1\,\mathrm{d}\mathcal{H}^{d-1} &\text{if } v=\chi_E \in BV(\mathbb{R}^d)\,, \|v\|_{L^1(\mathbb{R}^d)}=1\,,\\
        +\infty &\text{otherwise on }L^1(\mathbb{R}^d)\,.
        \end{cases}
    \end{align*}
    where $\partial^*E $ denotes the reduced boundary of the set $E$ and $\nu_E$ its unit outer normal.
\end{itemize}
}
\end{remark}

\subsection{The first eigenvalue of the combinatorial Laplacian}
In this section we introduce the first eigenvalue of the combinatorial Laplacian on subsets of $\mathbb{Z}^d$ and give some elementary properties.
Given $X\in \mathcal{X}_N$, we define its first eigenvalue as 
\begin{equation}\label{defin_discrete_F}
    \lambda_N\big(X\big) = \min\bigg\{  E_N(u)\colon  \ u(i)=0 \text{ in }\Z^d\setminus X\,,\quad \frac{1}{N} \sum_{i \in \mathbb{Z}^d} u^2(i)=1\bigg\}\,.
\end{equation}
If $X_N = N^{1/d}\Omega \cap \mathbb{Z}^d$ for some bounded and open $\Omega \subset \mathbb{R}^d$, by Remark~\ref{rem:energyscalings}{\rm (b)} and the Fundamental Theorem of $\Gamma$-convergence \cite{Bra2002,Dal1993}, we deduce that
\begin{align*}
\lim_{N \to +\infty} \lambda_N(X_N) = \lambda(\Omega)\,.
\end{align*}

\begin{definition} \label{def:minimalset} Given $N\in \mathbb{N}$ we say that $Y_N \in \mathcal{X}_N$ is a minimal set for $\lambda_N$ if 
\begin{align*}
\lambda_N(Y_N) \leq \lambda_N(X) \quad \text{for all } X\in \mathcal{X}_N\,.
\end{align*}
We moreover set
\begin{align}\label{def:m_lambda_N}
m_{\lambda,N} := {\inf}\{\lambda_N(X) \colon X \in \mathcal{X}_N\}\,.
\end{align}
\end{definition}

\subsection{Discrete rearrangements}

\noindent In this section we briefly recall some definitions and results related to discrete rearrangements and we refer the reader to \cite{HHH} and \cite{S-R} for more details. We denote by $D$ the following set of vectors
\begin{equation*}
    D=\bigg\{e_i, e_i+e_j, e_i-e_j \colon i,j=1,\dots,d \text{ and } i< j \bigg\}\,.
\end{equation*}
In what follows we enumerate the elements of $D$ as $\{v_1,\ldots,v_{d^2}\}$. For $k \in \mathbb{N}$ we set $v_k= ~{v_{((k-1) \,\mathrm{mod}\, d^2)+1}}$.

\begin{definition} \label{def:e-convex}
    Given a direction $e\in D$, we say that $X\subset\Z^d$ is $e$-convex if, for every  $x \in X$ and $K\in \mathbb{N}$ such that $x+Ke \in X$, we have that $x+ke \in X$ for all $k=1,\ldots,K-1$. We call a set $X\in \mathcal{X}$ direction-convex if it is $e$-convex for all $e \in D$.
\end{definition}

\noindent Given $e \in D$ we define
\begin{equation*}
\Pi_e =\begin{cases} \mathbb{Z}^d \cap \{x \colon \langle x,e\rangle=0\} &\text{if } e=e_k\,,\\
\mathbb{Z}^d \cap \big\{x \colon \langle x,e\rangle\in \{0,1\}\big\} &\text{if } e=e_i\pm e_j\,.
\end{cases}
\end{equation*}
Given $x_0 \in \mathbb{Z}^d$ we set $\Pi_e(x_0) = x_{0}+\Pi_{e}$.
Given $u \colon \mathbb{Z}^d \to [0,+\infty)$, $e \in D$, and $q \in \Pi_e$ we define $u^{q,e} \colon \mathbb{Z} \to [0,+\infty)$
\begin{align*}
u^{q,e}(t)=u({q+te}) \quad t\in \mathbb{Z}\,.
\end{align*}
Note that for every $e \in D$ each $i \in \mathbb{Z}^d$ can be uniquely written as $i=q+te$ for some $q\in \Pi_e$ and $t\in \mathbb{Z}$. Given $X \in \mathcal{X}$ we set $(X)^{q,e}=\{t\in \Z \colon q+te \in X\}$.



\begin{definition}\label{rearrangement_DEFIN} Let $u\colon \Z\to [0,+\infty)$ be a function with finite support. Let  $\{\alpha_i\}_{i\in \mathbb{N}}$ with  $\alpha_i \geq \alpha_{i+1}$ be the values taken by $u$. We define the symmetric decreasing rearrangement of $u$ as
        \begin{equation*}
            u^*(i)=\begin{cases}
                \alpha_{1-2i} &\text{if } i\leq 0\,,\\
                \alpha_{2i} &\text{if } i> 0\,.
            \end{cases}
        \end{equation*}
        Given $d\ge 2$, $e \in D$ and $u\colon \mathbb{Z}^d\to [0,+\infty)$ with finite support, we define the symmetric decreasing rearrangement of $u$ in direction $e$ as
        \begin{align*}
            u^{*e}({i}) := (u^{q,e})^*(t)\quad \text{for } i=q+te\,.
        \end{align*}
    Furthermore, setting $u^0=u$ and $u^{k}= (u^{k-1})^{*v_k}$, the symmetric decreasing rearrangement of $u$ is defined as
    \begin{align*}
        u^*:= \lim_{k\to +\infty} u^{k}\,.
    \end{align*}
    Given a finite set $X \in \mathcal{X}$ we define its symmetric rearrangement in direction $e$ as
    \begin{align*}
        R_e(X) := \mathrm{supp} (\chi_X^{*e})\,,
    \end{align*}
    while its symmetric rearrangement is defined as
    \begin{equation*}
 R(X):=  \mathrm{supp} (\chi_X^{*})\,.
    \end{equation*}
        We say that a finite set $X\in \mathcal{X}$ is symmetric if 
    \begin{align*}
    R_e(X) =X \quad \text{for all } e \in D\,.
    \end{align*}
\end{definition}    
\begin{remark}\label{finite_rearrangement_RMK} {\rm We remark the following properties of the Riesz rearrangement:
\begin{itemize}
    \item[{\rm (a)}]  In \cite{HHH} it has been shown that for $u \colon \mathbb{Z}^d \to [0,+\infty)$ with finite support $u^*$ is well-defined. Hence, if $X\subset \mathbb{Z}^d$ is a finite set $R(X)$ is well defined, too.
    \item [{\rm (b)}] By the definition of the rearrangement it follows that 
    \begin{align*}
        \sum_{i \in \mathbb{Z}^d} u^2(i) = \sum_{i \in \mathbb{Z}^d} (u(i)^*)^2\,.
    \end{align*}
    \item [{\rm (c)}] Note that for $X\in\mathcal{X}_N$ it holds that $\card R_e(X)=\card R(X)=N$.
\end{itemize}
}
\end{remark}

Before stating a discrete version of the Riesz rearrangement inequality, we introduce the notion of supermodular function.

\begin{definition}
    A function $G\colon \mathbb{R}\times\mathbb{R}\to \mathbb{R}$ is said to be supermodular if 
    \begin{equation*}
       G(x,y+t)+G(x+s,y)\leq  G(x+s,y+t)+G(x,y)   \quad \text{for any }x,y \in \mathbb{R},s,t>0\,.
    \end{equation*}
\end{definition}
The following discrete one-dimensional Riesz rearrangement inequality has been proved  in \cite[Proposition~4.3]{Haj1d}.

\begin{theorem}[1-dimensional Riesz rearrangement inequality]\label{thm:discrete_Riesz_1d}
    Let $u,v\colon\Z\to [0,+\infty)$ and let $H\colon \N\to [0,+\infty)$ be non-increasing. Let $G\colon\mathbb{R}\times\mathbb{R}\to \mathbb{R}$ be a supermodular function such that $G(0,0)=0$. Then,
   \begin{equation*}
         \sum_{i, j\in \Z} G\left(u(i),v(j)\right) H\left(|i-j|\right) \leq \sum_{i, j\in \Z} G\left(u^*(i),v^*(j)\right) H\left(|i-j|\right)\,.
    \end{equation*}
\end{theorem}

Similarly, the following theorem in higher dimensions has been proved in \cite[Theorem~1.2]{HHH}.

\begin{theorem}[Riesz rearrangement inequality]\label{thm:discrete_Riesz_rearrangement}
    Let $u,v\colon\Z^d\to [0,+\infty)$ be two functions with finite support and let $H\colon \N\to [0,+\infty)$ be non-increasing. Let $G\colon [0,+\infty)\times[0,+\infty)\to [0,+\infty)$ be a supermodular function such that $G(0,0)=0$. Then,
   \begin{equation*}
         \sum_{i, j\in \Z^d} G\left(u(i),v(j)\right) H\left(\norm{i-j}_{1}\right) \leq \sum_{i, j\in \Z^d} G\left(u^*(i),v^*(j)\right) H\left(\norm{i-j}_{1}\right)\,.
    \end{equation*}
\end{theorem}

The following corollary is a consequence of the previous two Theorems.

\begin{corollary} \label{cor:sym}{\rm  Let $u \colon \mathbb{Z}^d\to [0,+\infty)$ have finite support.
\begin{itemize}
\item[{\rm (a)}]  For any $e\in D$, the symmetric decreasing rearrangement of $u$ in direction $e$ satisfies
\begin{align*}
    &\sum_{(i,j)\in \mathcal{N}} |u^{*e}(i) - u^{*e}(j)|^2 \leq  \sum_{(i,j)\in \mathcal{N}}|u(i) - u(j)|^2\,.
\end{align*}
In particular, for any $X \in\mathcal{X}_N$ it follows that $\lambda_N\big(R_e(X)\big) \leq \lambda_N(X)$. 
\item[{\rm (b)}] The symmetric rearrangement of $u$ satisfies
\begin{align*}
    &\sum_{(i,j)\in \mathcal{N}} |u^*(i)-u^*(j)|^2 \leq \sum_{(i,j)\in \mathcal{N}} |u(i)-u(j)|^2\,.
\end{align*}
In particular, for any $X\in \mathcal{X}_N$, it follows that $\lambda_N\big(R(X)\big) \leq \lambda_N(X)$.
\end{itemize}
}
\end{corollary}

\begin{proof}  Let $u\colon \mathbb{Z}^d \to [0,+\infty)$ with finite support. First of all, note that the function $G\colon \mathbb{R} \times \mathbb{R} \to \mathbb{R}$ defined by $G(x,y)=-|x-y|^2$ is supermodular. \\
\noindent \emph{Proof of {\rm (a)}.}  Let $e \in D$. In this proof we distinguish different cases and make use of Theorem~\ref{thm:discrete_Riesz_1d} with $G(x,y)=-|x-y|^2$ and different choices for $H$, $u$, and $v$. 

If $e = e_n$ for some $n=1,\ldots,d$, the inequality follows choosing $H(t)=\chi_{\{0,1\}}(t)$, and $u=v=u^{q,e}(t)$ with $t \in \mathbb{Z}$ and $q \in \Pi_e$ for the interactions in direction $e_n$, while for the interactions in direction $e_k$, $k\neq n$, one chooses $H(t)=\chi_{\{0\}}(t)$, and $u=u^{q,e}(t)$, $v=u^{q+e_k,e}(t)$, $k\neq n$ with $t \in \mathbb{Z}$ and $q \in \Pi_e$.

 In the case $e=e_n+e_k$, the inequality for the direction $e_j$, $j \notin \{k,n\}$ follows as above  by choosing $H(t)=\chi_{\{0\}}(t)$, $u=u^{q,e}(t)$ and $v=u^{q+e_j,e}(t)$, with $t \in \mathbb{Z}$ and $q \in \Pi_e$. For the directions $e \in \{e_n, e_k\}$ the result follows from \cite[Section~5.2]{S-R}  by noting that the lines in direction $e=e_n+e_k$ are contained in the $2$-dimensional plane $\mathrm{span}_\mathbb{Z}\{e_n,e_k\}$.\\
\noindent \emph{Proof of {\rm (b)}.}  This follows from Theorem~\ref{thm:discrete_Riesz_rearrangement} with $G(x,y)=-|x-y|^2$, $H(t)=\chi_{\{0,1\}}(t)$, and $u=v$.
This concludes the proof of the corollary.
\end{proof}

\section{Properties of optimal sets} \label{sec:prop_optimal_sets}
\noindent In this section we prove some properties of functions minimizing \eqref{defin_discrete_F} and of discrete sets optimizing \eqref{def:m_lambda_N}. We start by introducing the notion of connectedness of subsets of $\mathbb{Z}^d$.
 The following elementary lemma is proved here for the reader's convenience.

\begin{lemma} \label{lem:boundedness-of-lambdaN} There exists $C_d>0$ such that for all $N\in \mathbb{N}$ there holds
\begin{align*}
m_{\lambda,N}
 \leq C_d\,.
\end{align*}
\end{lemma}
\begin{proof} For each $N\in \mathbb{N}$, we construct a competitor $X \in \mathcal{X}_N$ such that
\begin{align} \label{ineq:lambdaN}
 \lambda_N(X) \leq C_d
\end{align}
for some $C_d>0$ independent of $N \in \mathbb{N}$. We only need to prove \eqref{ineq:lambdaN} for $N$ large enough. First, we perform the construction for $N = (2k+1)^d$ for $k\in \mathbb{N}$. To this end, let $X = [-k,k]^d \cap \mathbb{Z}^d$ and $u \colon \mathbb{Z}^d \to \mathbb{R}$ be defined by
\begin{align*}
u(i) = C_{d,N} (k-l) \quad \text{ for all } i \in \partial [-l,l]^d \cap \mathbb{Z}^d\,, 0\leq l \leq k\,,
\end{align*}
where $C_{d,N} >0$ is such that 
$
\sum_{i \in \mathbb{Z}^d} u^2(i) =N
$. We show that there exists $c_d >0$ such that
\begin{align}\label{ineq:CdN}
\frac{1}{c_d} N^{-\frac{1}{d}}\leq  C_{d,N} \leq c_d N^{-\frac{1}{d}}\,.
\end{align}
 Note that $  l^{d-1} \leq \card\left(\partial [-l,l]^d \cap \mathbb{Z}^d\right) \leq 2d  l^{d-1} $ and therefore for $k$ large enough an elementary calculation shows
\begin{align*}
8^{-d} C_{d,N}^2k^{d+2}   \leq C_{d,N}^2 \sum_{l=0}^{k} (k-l)^2l^{d-1}\leq \sum_{i \in \mathbb{Z}^d} u^2(i) \leq  C_{d,N}^2 2d\sum_{l=0}^k (k-l)^2l^{d-1} \leq 2d C_{d,N}^2 k^{d+2}\,.
\end{align*}
Recalling that $\sum_{i \in \mathbb{Z}^d} u^2(i) =N$ the inequalities in \eqref{ineq:CdN} follow. Finally, using \eqref{ineq:CdN}, we obtain
\begin{align*}
E_N(u) = N^{-\frac{d-2}{d}}\sum_{(i,j) \in \mathcal{N}}|u(i)-u(j)|^2 \leq N^{-\frac{d-2}{d}} C_{d,N}^2 2d\sum_{l=0}^k  l^{d-1} \leq 2d c_d^2 N^{-1} k^d \leq  2d c_d^2 =:C_d  \,.
\end{align*}
Setting $ N_{\mathrm{min}}= (2k+1)^d \leq N < (2k+3)^d-1=N_{\mathrm{max}}$ the claim follows by taking as a competitor the one constructed before for $N=N_{\mathrm{min}}$, rescaling it to keep the mass constraint and eventually noting that  $ N_{\mathrm{max}} \leq  2^d N_{\mathrm{min}}$.
\end{proof}

\begin{definition} \label{def:connectedness}
    We say that $X\subset\Z^d$ is {connected} if, for all $p,q \in X$, there exist $M\in\N$ and points $p_{1},\dots,p_{M}\in X$ such that
        \begin{equation*}
            \begin{cases}
                p_0=p, \quad p_{M}=q\,,\\
                (p_k,p_{k-1}) \in \mathcal{N} \quad \forall k=1,\dots, M\,.
            \end{cases}           
        \end{equation*}
        Given $X \in \mathcal{X}$ we call a connected component of $X$ any maximal (with respect to set inclusion) connected subset of $X$.  
\end{definition}

\begin{proposition} \label{analytical_properties_u_PROP}
    Let $X\subset\Z^d$ be a connected set. Then, the following holds:
    \begin{enumerate}
        \item[{\rm (a)}] any minimizer $u$ of (\ref{defin_discrete_F}) is such that $u(i)> 0$ for each $i\in X$;
        \item[{\rm (b)}] there exists a unique function that minimizes (\ref{defin_discrete_F}).
    \end{enumerate}
\end{proposition}

\begin{proof}[Proof of {\rm (a)}.] The fact that $u(i)$ is of constant sign follows directly from
\begin{align*}
    \big|u(i)-u(j)\big|\geq\big||u(i)|-|u(j)|\big| \quad \text{for all } u(i),u(j) \in \mathbb{R}\,.
\end{align*}
Without loss of generality we can assume that $u(i)\geq 0$ on $X$. Let us suppose by contradiction that there exists a minimizer $u$ that vanishes on $A \subseteq X$ (maximal with respect to set inclusion) such that $M:=\# A<\# X$. Then, for $t\in(0,1)$ we consider $u^t \colon \mathbb{Z}^d \to \mathbb{R}$ defined by
\begin{align*}
u^t(i):=\begin{cases}
	tu(i) &\text{if } i \in \mathbb{Z}^d\setminus A\,,\\
	\sqrt{\frac{N}{M}(1-t^2)} &\text{if } i\in A\,.
        \end{cases}
\end{align*}
Let us point out that, for each $t\in (0,1)$, it holds that
\begin{align*}
\frac{1}{N}\sum_{i\in \mathbb{Z}^d}(u^t(i))^2= \frac{1}{N}\sum_{i\in \mathbb{Z}^d\setminus A} t^2 u^2(i) + \frac{1}{N}\sum_{i\in A} \bigg(\sqrt{\frac{N}{M}(1-t^2)}\bigg)^2 = t^2 + (1-t^2) =1\,.
\end{align*}
Hence, $u^t$ is a competitor for the minimum in (\ref{defin_discrete_F}). However, the Dirichlet energy of $u^t$ is
\begin{align*}
    \sum_{(i,j)\in \mathcal{N}} |u^t(i)-u^t(j)|^2  &= \underset{i,j \in \mathbb{Z}^d\setminus A}{\sum_{(i,j) \in \mathcal{N} }} |u^t(i)-u^t(j)|^2 +2\underset{i \in \mathbb{Z}^d\setminus A, j \in A}{\sum_{(i,j) \in \mathcal{N} }}  |u^t(i)-u^t(j)|^2+ \underset{i,j\in A}{\sum_{(i,j) \in \mathcal{N} }}  |u^t(i)-u^t(j)|^2\\
    &= t^2\underset{i,j \in \mathbb{Z}^d\setminus A}{\sum_{(i,j) \in \mathcal{N} }} |u(i)-u(j)|^2 + 2\underset{i \in \mathbb{Z}^d\setminus A, j \in A}{\sum_{(i,j) \in \mathcal{N} }}\left|t\,u(i) -\sqrt{\frac{N}{M}(1-t^2)}\right|^2\\
    &=t^2 D(u) +2\sqrt{\frac{N}{M}(1-t^2)}\underset{i \in \mathbb{Z}^d\setminus A, j \in A}{\sum_{(i,j) \in \mathcal{N} }}\left(-2t\, u(i)+\sqrt{\frac{N}{M}(1-t^2)}\right)\,.
\end{align*}
Note that, since $X$ is connected and $A$ is maximal with respect to set inclusion, there exist $(i,j) \in \mathcal{N}$ with $i \in \mathbb{Z} \setminus A$ and $j \in A$ such that $u(i_0)>0$. Since the sum (being finite) on the right-hand side of the previous equation is continuous in $t$, and for $t=1$ we have 
\begin{equation*}
  \underset{i \in \mathbb{Z}^d\setminus A, j \in A}{\sum_{(i,j) \in \mathcal{N} }} -2 u(i) \leq -2u(i_0) <0\,,
\end{equation*}
there exists $t \in (0,1)$ such that 
\begin{align*}
2\sqrt{\frac{N}{M}(1-t^2)}\underset{i \in \mathbb{Z}^d\setminus A, j \in A}{\sum_{(i,j) \in \mathcal{N} }}\left(-2t\, u(i)+\sqrt{\frac{N}{M}(1-t^2)}\right) <0\,.
\end{align*}
This yields a contradiction to the minimality of $u$. \\

\noindent \emph{Proof of {\rm (b)}.} Let us suppose that $u,v$ are both minimizers of $\lambda_{N}(X)$ in \eqref{defin_discrete_F}. Thanks to part {\rm (a)} we can assume that $u,v>0$ on $X$.  Given  $t\in(0,1)$, let us define $w^t\colon \mathbb{Z}^d \to \mathbb{R}$ as
\begin{align*}
w^t(i):=(tu^2(i)+(1-t)v^2(i))^{1/2}\,.
\end{align*}
 Since $\frac{1}{N}\sum_{i\in X}|w^t_i|^2=1$, $w^t$ is a competitor for $\lambda_N(X)$. Moreover, for $(i,j) \in \mathcal{N}$ it holds that
    \begin{align*}
        |w^t(i)-w^t(j)|^2 &= \big(t\,u^2(i)+(1-t)v^2(i)\big)+\big(t\,u^2(j)+(1-t)v^2(j)\big)\\
        &\qquad -2(t^2u^2(i) u^2(j) +(1-t)^2 v^2(i) v^2(j) +t(1-t)(u^2(i) v^2(j) +u^2(j) v^2(i)))^{1/2}\\
        &= t\big(u^2(i)+u^2(j)\big)+(1-t)\big(v^2(i)+v^2(j)\big)\\
        &\qquad-2(t^2u^2(i) u^2(j) +(1-t)^2v^2(i) v^2(j) +t(1-t)(u^2(i) v^2(j) +u^2(j) v^2(i)))^{1/2}\,.
    \end{align*}
From the previous calculation we deduce that
    \begin{align*}
        \begin{split}
        0 &\leq E_N(w^t) - \lambda_{N}(X)=E_N(w^t)- t E_N(u) -(1-t)E_N(v)  \\
        &= \sum_{(i,j) \in \mathcal{N}} |w^t(i)-w^t(j)|^2 \\
        &\qquad\qquad\qquad- \sum_{(i,j) \in \mathcal{N}} \bigg(t(u^2(i)+u^2(j))+(1-t)(v^2(i)+v^2(j)) -2tu(i) u(j) -2(1-t)v(i) v(j)\bigg) \\
        &= 2\sum_{(i,j) \in \mathcal{N}}  \bigg(tu(i) u(j) +(1-t)v(i) v(j) \\
        &\qquad\qquad\qquad-(t^2u^2(i) u^2(j) +(1-t)^2v^2(i) v^2(j) +t(1-t)(u^2(i) v^2(j) +u^2(j) v^2(i)))^{1/2}  \bigg)\,.
	\end{split}
    \end{align*}
Note that each term in the last sum is non-positive and equals zero if and only if $u_{i}v_{j} = u_{j}v_{i}$, as an elementary computation shows. Hence, we have that $\frac{u(i)}{v(i)}=\frac{u(j)}{v(j)}$. Since $X$ is connected, $u(i)=\alpha v(i)$ for some $\alpha>0$ and for all $i\in X$. By the constraint on the norm, we infer that $\alpha=1$ which concludes the proof of (b).
\end{proof}

The following geometric Lemma will be used in the proof of Proposition~\ref{discrete_minimum_properties_PROP}.

\begin{lemma} \label{lemma:translation} Let $X,Y \subset \mathbb{Z}^d$ have finite cardinality. Then, there exists $\tau \in \mathbb{Z}^d$ such that there exists a unique pair $(x,y) \in \mathcal{N}$ with $x \in X$ and $y \in Y+\tau$. 
\end{lemma}

\begin{proof} The proof proceeds in two steps. First we show that for any $X,Y \in \mathbb{Z}^d$ with finite cardinality there exists a translation $\sigma \in \mathbb{Z}^d$ such that
\begin{align} \label{eq:translation-intersection}
\begin{split}
 {\rm (i)}& \, X \cap (Y+\sigma)=\{z\} \text{ for some } z \in \mathbb{Z}^d\,; \\   {\rm (ii)}& \, X \subset \{x \in \mathbb{Z}^d\colon x_d\geq z_d\} \quad \text{ and } \quad Y+\sigma \subset \{y \in \mathbb{Z}^d\colon y_d\leq z_d\}\,.
\end{split}
\end{align}
In order to prove the claim of the Lemma,  assuming \eqref{eq:translation-intersection}, it suffices to choose $\tau=\sigma-e_d$. We are left with the proof of \eqref{eq:translation-intersection}. To this end, we proceed by induction on the dimension. For $d=1$ the proof is trivial. We assume now that the statement is proven for $d-1$ and proceed to prove it for $d$. To this end, we set
\begin{align*}
m_d(X) := \min_{x \in X} x_d \quad \text{ and } \quad M_d(Y) := \max_{y \in Y} y_d\,
\end{align*}
and define $\tau_d := -M_d(Y) +m_d(X)$. We apply the induction assumption for $X_{d-1}=X \cap (\mathbb{Z}^{d-1}\times\{m_d(X)\} )$ and $Y_{d-1}=Y\cap (\mathbb{Z}^{d-1}\times\{M_d(Y)\} )$ to find $\hat{\tau} \in \mathbb{Z}^{d-1}$ such that \eqref{eq:translation-intersection} holds true (with a slight abuse of notation we identify $X_{d-1}$ and $Y_{d-1}$ as subsets of $\mathbb{Z}^{d-1}$). Eventually, the translation $\tau=(\hat{\tau},\tau_d)$ satisfies \eqref{eq:translation-intersection}.
\end{proof}

\begin{proposition}\label{discrete_minimum_properties_PROP} Let $Y_N\in \mathcal{X}_N$ be a minimal set for $\lambda_N$ according to Definition~\ref{def:minimalset}. Then, the following properties hold:
    \begin{itemize}
        \item[{\rm (a)}] $Y_N$ is connected;
        \item[{\rm (b)}] $Y_N$ is direction-convex; 
        \item[{\rm (c)}]  There exists a constant $C_d >0$ (independent of $N\in \mathbb{N}$) such that if $d=2$ or if $d>2$ and $Y_N$ is symmetric according to Definition~\ref{rearrangement_DEFIN}, then
        \begin{align*}
        \mathrm{diam}(Y_N) \leq C_d N^{1/d}\,.
        \end{align*} 
   \end{itemize}
\end{proposition}
\begin{remark} \label{rem:symmetric} {\rm We note that, thanks to \eqref{cor:sym}({\rm b}), for any $d\geq 2$ and $N \in \mathbb{N}$ there exists a symmetric minimal set $X_N \in \mathcal{X}_N$.}
\end{remark}
\begin{proof}[Proof of Proposition~\ref{discrete_minimum_properties_PROP}]
For the proof of {\rm (a)} we argue by contradiction assuming that $Y_N$ is not connected. To fix the ideas (the general case follows by the same argument) we assume that $Y_N = C^1\cup C^2$ with $C^1 ,C^2 \subset \mathbb{Z}^d$ two disjoint connected components of $Y_N$ according to Definition \ref{def:connectedness}. We denote by $u$ a non-negative function realizing the minimum for $\lambda_N(Y_N)$ according to definition \eqref{defin_discrete_F}, that is
\begin{align*}
E_N(u) = \lambda_N(Y_N)\text{ and } u =0 \text{ on } \mathbb{Z}^d\setminus Y_N\,.
\end{align*}
Let $\tau \in \mathbb{Z}^d\setminus \{0\}$ be the translation given in Lemma~\ref{lemma:translation} for $X=C^1$ and $Y=C^2$ such that there exists a unique pair $(x,y) \in \mathcal{N}$ with $x \in C^1$ and $y \in C^2+\tau$. The set $\hat{Y}_N= C^1\cup (C^2+\tau)$ is connected and such that $\#\hat{Y}_N=N$. We define 
\begin{align*}
\hat{u}(i) = \begin{cases} u(i) &\text{if } i \notin C^2+\tau\,,\\
u({i-\tau})  &\text{if } i \in C^2+\tau\,,
\end{cases}
\end{align*}
and observe that, since $\hat{u}$ is a competitor for $\lambda_N(\hat{Y}_N)$, 
\begin{align*}
E_N(\hat{u}) \geq \lambda_N(\hat{Y}_N)\,.
\end{align*}
Furthermore, as $C^1$ and $C^2$ are two disjoint connected components of $Y_N$, we have
\begin{align*}
E_N(\hat{u}) &= E_{N}(\hat{u}\chi_{C^{1}}) + E_{N}(\hat{u}\chi_{C^{2}+\tau}) + |\hat{u}({x})-\hat{u}({y})|^{2} - |\hat{u}({x})|^{2} - |\hat{u}({y})|^{2}\\
&= E_{N}(u) + |u({x})-u({y-\tau})|^{2} - |u({x})|^{2} - |u({y-\tau})|^{2}\\
&\le E_{N}(u)\,,
\end{align*}
with equality if and only if $\min\left(u({x}),u({y-\tau})\right)=0$. This shows that $\hat{u}$ is a minimizer and, according to Proposition \ref{analytical_properties_u_PROP}, it satisfies in particular $u(x),u({y-\tau})>0$, which gives a contradiction to the minimality of $Y_N$.\\

\noindent\emph{Proof of {\rm(b)}.} This is a consequence of Corollary~\ref{cor:sym}~{\rm (a)} and Remark~\ref{finite_rearrangement_RMK}~{\rm (c)}.\\

 \noindent\emph{Proof of {\rm(c)}.} The case $d=2$ follows from \cite[Proposition~6.6]{S-R}. It remains to prove the case $d>2$ for $Y_N$ a symmetric and minimal set for $\lambda_N$. We claim that there exists $C_d>0$ such that for all $i \in Y_N$ there holds
\begin{align}\label{ineq:coordinate_estimate}
|i_n| \leq C_d N^{\frac{1}{d}} \text{ for all } n=1,\ldots,d\,.
\end{align}
Clearly, the claim implies ({\rm c}). We are thus left to prove \eqref{ineq:coordinate_estimate}. Assume that there exists $i \in Y_N$ such that (without loss of generality) $|i_d| \geq \kappa$. We then show that there exists $C_d>0$ such that
\begin{align}\label{ineq:YNestimate}
\#Y_N \geq C_d \kappa^d\,.
\end{align}
 Since $R_{e_d}(Y_N) =Y_N$, there exist $i^* \in Y_N$ such that 
\begin{align*}
 \begin{cases} i^*_n =i_n &\text{if } n \in \{1,\ldots,d-1\}\,,\\
|i^*_d-i_d| \leq 1 \,.
\end{cases}
\end{align*}
Thanks to part ({\rm b}) $Y_N$ is direction-convex and therefore contains the segment $[i,i^*] \cap \mathbb{Z}^d$. We can assume that this segment contains the origin. In fact, if it is not the case, there exists an index $k \in \{1,\ldots,d-1\}$ such that for all $j \in [i,i^*] \cap \mathbb{Z}^d$ we have $j_k \neq 0$. At this point, we also have
\begin{align*}
R_{e_k} ([i,i^*] \cap \mathbb{Z}^d)=\{R_{e_k}j \colon j \in ([i,i^*] \cap \mathbb{Z}^d)\} \subset Y_N\,,
\end{align*}
where 
\begin{align*}
R_{e_k}j_n = \begin{cases} j_n &\text{if } n \neq k\,,\\
-j_k &\text{if } n=k\,.
\end{cases}
\end{align*}
Exploiting the direction-convexity of $Y_N$ for the direction $e_k$, we conclude that 
$
([i,i^*] -i_k)\cap \mathbb{Z}^d \subset Y_N
$ and passes through the origin. This implies that $[i,i^*] \cap \mathbb{Z}^d\supset \{ n e_d \colon n\in \mathbb{N}\,, |n|\leq \kappa-1\}$. Again by symmetry with respect to $e=e_j+e_d$ for $j\neq d$, we have that $\{ n e_j \colon n\in \mathbb{N}\,, |n|\leq \kappa-1\} \subset Y_N$ for all $j=1,\ldots,d-1$. Due to the directional convexity with respect to $e=e_j+e_k$ for all $j\neq  k$, we conclude that 
\begin{align*}
C:=\mathrm{conv}\left(\bigcup_{j=1}^d \{ n e_j \colon n\in \mathbb{N}\,, |n|\leq \kappa-1\}  \right) \cap \mathbb{Z}^d \subset Y_N\,.
\end{align*}
Since $C$ satisfies $\# C= C_d \kappa^d$ we obtain \eqref{ineq:YNestimate}. Since $\#Y_N=N$, this implies \eqref{ineq:coordinate_estimate}. 
\end{proof}

\section{Maximal Fluctuation Estimates} \label{sec:fluctuations}

In this section we prove the main result of the paper, that is stated in the following theorem.

\begin{theorem}\label{main_THM} Let $\{\alpha_N\}_N\subseteq (0,+\infty)$ be such that $\sup_{N} \alpha_N <+\infty$ and let $X \in \mathcal{X}_N$ satisfy
\begin{align}\label{ineq:Xalmostmin}
\lambda_N(X) \leq m_{\lambda,N} + \alpha_N\,.
\end{align}
Then, there exists $C_d>0$ such that
\begin{align}\label{ineq:fluctuation}
\begin{split}
\min_{z \in \Z^d}\card\left(X \Delta ((B_{r_N} \cap \mathbb{Z}^d)+z)\right) &\leq C_d N\left(
   N^{-\frac{1}{2d}} P_N(X)^{\frac{1}{2}}  + \alpha_N^{\frac{1}{2}} +N^{-\frac{1}{d}} \right)\,,
\end{split}
\end{align}
where  $r_N>0$ is such that $|B_{r_N}|=N$. In particular, if   $X_N \in \mathcal{X}_N$ is such $\sup_{N}\mathrm{P}_N(X_N)<+\infty$ and 
\begin{align*}
\lambda_N(X_N) \leq m_{\lambda,N} + \alpha_N\,,
\end{align*}
there holds
\begin{equation}\label{ineq:quasiminimi}
\min_{z \in \mathbb{Z}^d}\card \left(X_N \Delta ((B_{r_N} \cap \mathbb{Z}^d)+z)) \right) \leq C_d N \left( N^{-\frac{1}{2d}} + \alpha_N^{\frac{1}{2}}\right) \,.
\end{equation}
\end{theorem}

\begin{remark}\label{rem:Perimter-estimate}  {\rm We observe that the estimate
\begin{align*}
\mathrm{sup}_N \mathrm{P}_N(X_N) <+\infty
\end{align*} 
holds true for any $X_N \in \mathcal{X}_N$ minimal set, if $d=2$ or if  $d>2$ and $X_N$ is symmetric. In fact, the perimeter estimate is implied by Proposition~\ref{discrete_minimum_properties_PROP}({\rm b}),({\rm c}) as for directional-convex sets there holds
\begin{align*}
c_1\mathrm{diam}(X_N)^{d-1}\leq \mathrm{P}(X_N) \leq c_2\mathrm{diam}(X_N)^{d-1}\,
\end{align*}
for some dimensional constants $c_1,c_2>0$. Thanks to Remark~\ref{rem:symmetric}, for any $d\geq 2$ and $N\in \mathbb{N}$ there exists a minimal set   $X_N \in \mathcal{X}_N$ that satisfies the above estimate. }
\end{remark}

\subsection{Embedding of the discrete problem into the continuum setting} \label{sec:zetaandaffine}

 In order to follow the strategy of the proof of the main result highlighted in the Introduction, here we show how to embed the discrete problem in a continuum setting. This amounts to properly extending both sets $X \in \mathcal{X}$ and functions $u\colon \mathbb{Z}^d\to \mathbb{R}$. To this aim, we exploit the Kuhn decomposition of a cube introduced in Section~\ref{subsec:Kuhn}.  We begin by extending sets and introducing the map $\zeta \colon \mathcal{X} \to \mathfrak{M}(\mathbb{R}^d)$, where $\mathfrak{M}(\mathbb{R}^d)$ denotes the measurable subsets of $\mathbb{R}^d$.  Given $X \in \mathcal{X}$ we set
\begin{align}\label{def:zeta}
\zeta(X) = \bigcup_{i \in X} \bigcup_{T \in \mathcal{T}(i)} T\,.
\end{align}
\begin{definition} \label{def:piecewiseaffine}
Given a function $u \colon \mathbb{Z}^d \to \mathbb{R}$ we introduce the function $\hat{u}\colon\mathbb{R}^d\to\mathbb{R}$ as the function who is affine on each Kuhn simplex $T \in \mathcal{T}$ and interpolates the values of $u$ on the points of the lattice.
\end{definition}
 
The next lemma estimates the error in the measure one makes when passing from a discrete set $X$ to its continuum representation $\zeta(X)$.

\begin{lemma}\label{lemma_measure_zeta}
Let $N\in\mathbb{N}$, $X \in \mathcal{X}_N$ and let $\zeta$ be the map defined in \eqref{def:zeta}. Then, there exists $C_d>0$ such that
\begin{equation}\label{measure_zeta_estimate}
    N\leq \big|\zeta(X)\big| \leq N+C_d N^{\frac{d-1}{d}} P_N(X)\,.
\end{equation}
 In particular, $\big|\zeta(X)\big| \leq C_d N$. 
\end{lemma}
\begin{proof}
    The lower bound is straightforward since, by the definition of $T \in \mathcal{T}(i)$, we have
    \begin{equation*}
        Q_1(i)\subseteq  \bigcup_{T \in \mathcal{T}(i)} T\,.
    \end{equation*}
Given $X \in \mathcal{X}_N$, we first notice that
\begin{align*}
    \zeta(X) \subseteq \bigcup_{i \in (X)_{\sqrt{d}}} Q_1(i)\,
\end{align*}  
(recall that $(X)_{\sqrt{d}} = \{ i\in \mathbb{Z}^d \colon  \mathrm{dist}(i,X) \leq \sqrt{d}\}$).
As a consequence, noting that the union on the right hand side is a disjoint union, we obtain
\begin{align*}
|\zeta(X)| \leq \sum_{i \in (X)_{\sqrt{d}}} |Q_1(i)| =\card X + \card((X)_{\sqrt{d}}\setminus X)=N +\card((X)_{\sqrt{d}}\setminus X) \,.
\end{align*}
Next, note that $i \in (X)_{\sqrt{d}}\setminus X$ only if there exists $j \in X$ with $\mathrm{dist}(i,j) \leq \sqrt{d}$ and $\mathrm{val}(j)\geq 1$. Moreover, for each $j\in \Z^{d}$ there exists a constant $C_d>0$ (independent of $j$) such that
\begin{align*}
\card\{i \in \mathbb{Z}^d\setminus X \colon \mathrm{dist}(i,j) \leq \sqrt{d}\} \leq C_d\,.
\end{align*}
Hence,
\begin{align*}
\card\left( (X)_{\sqrt{d}}\setminus X \right) \leq  C_d\sum_{x \in X} \mathrm{val}(x) = C_d N^{\frac{d-1}{d}} P_N(X)\,,
\end{align*}
which concludes the proof of \eqref{measure_zeta_estimate}.  The last part of the statement follows by construction. 
\end{proof}

 In what follows, we compare the energy functionals needed to define the discrete and continuum eigenvalues. 

\begin{lemma}\label{utilde_properties_LEMMA} Let $u \colon \mathbb{Z}^d \to \mathbb{R}$ be such that $\mathrm{supp}(u)=X\in\mathcal{X}_N$. The function $\hat{u}$ in Definition~\ref{def:piecewiseaffine} has the following two properties:
    \begin{enumerate}
        \item[{\rm (a)}]  \begin{equation*}
             \int_{\Rd} |\nabla \hat{u}(x)|^2\dx = N^{\frac{d-2}{d}}E_N(u)\,,
            \end{equation*}
        \item[{\rm (b)}] 
        \begin{equation*}
        \begin{split}
        \int_{\Rd} |\hat{u}(x)|^2\dx &\geq  \sum_{i \in \mathbb{Z}^d} |u(i)|^2 - 2\sqrt{d} \left(\sum_{i \in \mathbb{Z}^d} |u(i)|^2\right)^{1/2} N^{\frac{d-2}{2d}}E_N(u)^{\frac{1}{2}}\,.
        \end{split}
        \end{equation*}
    \end{enumerate}
\end{lemma}

\begin{proof} We first prove {\rm (a)}. To this end, we prove for all $k \in \{1,\ldots,d\}$ that
\begin{align} \label{eq:partialktilde}
\sum_{i \in \mathbb{Z}^d} |u({i+e_k}) - u(i)|^2 = \int_{\mathbb{R}^d} |\partial_k \hat{u}|^2\,\mathrm{d}x\,.
\end{align}
Note that once \eqref{eq:partialktilde} is proven, {\rm (a)} follows by summing over $k \in \{1,\ldots,d\}$. Let us prove \eqref{eq:partialktilde}. To this end, fix $k \in \{1,\ldots,d\}$ and let $i \in \mathbb{Z}^d$. We denote by $\mathcal{T}_{i,k}$ the set of simplices that have $[i,i+e_k]$ as an edge and observe that, thanks to Proposition~\ref{Kuhn_properties_PROP}(c),  $\#\mathcal{T}_{i,k}=d!$. For $T \in \mathcal{T}_{i,k}$ we note that 
\begin{align*}
\partial_k \hat{u}(x) = u({i+e_k})-u(i) \text{ for all } x \in T\,.
\end{align*}
Hence, by Proposition~\ref{Kuhn_properties_PROP} and the fact that $\mathring{T}_1 \cap \mathring{T}_2=\emptyset$ for all $T_1, T_2 \in \mathcal{T}$, $T_1 \neq T_2$,  we obtain
\begin{align*}
\int_{\bigcup T \in \mathcal{T}_{i,k}} |\partial_k \hat{u}(x) |^2\,\mathrm{d}x &=\sum_{T \in \mathcal{T}_{i,k}}\int_{T} |\partial_k \hat{u}(x) |^2\,\mathrm{d}x =   \sum_{T \in \mathcal{T}_{i,k}}|T|\cdot |u({i+e_k})-u(i)|^2\\
&= \frac{1}{d!} \#\mathcal{T}_{i,k} |u({i+e_k})-u(i)|^2 = |u({i+e_k})-u(i)|^2 \,.
\end{align*}
Summing over $i \in \mathbb{Z}^d$ we eventually get
\begin{align*}
\sum_{i \in \mathbb{Z}^d}|u({i+e_k})-u(i)|^2=\sum_{i \in \mathbb{Z}^d} \int_{\bigcup T \in \mathcal{T}_{i,k}} |\partial_k \hat{u}(x) |^2\,\mathrm{d}x=\int_{\mathbb{R}^d} |\partial_k \hat{u}(x) |^2\,\mathrm{d}x\,,
\end{align*}
which concludes the proof of {\rm (a)}. \\ \noindent Let us prove {\rm (b)}. Let $T \in \mathcal{T}$, $x \in T$ and let $z \in T$ be a vertex of $T$.  By definition of $\hat{u}$ we have
\begin{align*}
|\hat{u}(x)-u(z)|=|\hat{u}(x)-\hat{u}(z)| \leq \sqrt{d} \, |\nabla \hat{u}_{|_T}| \,.
\end{align*}
As a result, by the triangular inequality, we have
\begin{align*}
    u^{2}(x) \geq \big| |\hat{u}(x)-u(z)|-|u(z)| \big|^2 \geq u^{2}(z) - 2\sqrt{d}|\nabla \hat{u}_{|_T}| |u(z)|\,.
\end{align*}
Again by the very definition of $\hat{u}$ we observe that 
\begin{align*}
|\nabla \hat{u}_{|_T}|^2 \leq \underset{i,j\in \mathbb{Z}^d \cap T}{\sum_{(i,j) \in \mathcal{N}}}|u(i)-u(j)|^2\,.
\end{align*}
Making use of the previous inequality, Proposition~\ref{Kuhn_properties_PROP}(a) and Cauchy-Schwarz inequality, we have
\begin{align}\label{lemma:poincare}
\nonumber\sum_{i \in \mathbb{Z}^d} \sum_{\pi \in \mathcal{P}_d} \int_{T_\pi(i)} \hat{u}^2(x)\,\mathrm{d}x &\geq \sum_{i \in \mathbb{Z}^d} \sum_{\pi \in \mathcal{P}_d} |T_\pi(i)|\left( u^{2}(i) - 2\sqrt{d}|\nabla \hat{u}_{|_{T_\pi(i)}}| |u(i)|\right) \\&\nonumber\geq \sum_{i \in \mathbb{Z}^d} u^2(i)-   \sum_{i \in \mathbb{Z}^d} \sum_{\pi \in \mathcal{P}_d} \frac{1}{d!} 2\sqrt{d} \, |u(i)| \Big( \underset{h,k\in \mathbb{Z}^d \cap T_\pi(i)}{\sum_{(h,k) \in \mathcal{N}}}|u(h)-u(k)|^2\Big)^{\frac{1}{2}} \\&\geq \sum_{i \in \mathbb{Z}^d} u^2(i)- \frac{2\sqrt{d}}{d!} \Big(\sum_{i \in \mathbb{Z}^d} \sum_{\pi \in \mathcal{P}_d}u^2(i) \Big)^{\frac{1}{2}}  \Big(\sum_{i \in \mathbb{Z}^d} \sum_{\pi \in \mathcal{P}_d} \underset{h,k\in \mathbb{Z}^d \cap T_\pi(i)}{\sum_{(h,k) \in \mathcal{N}}}|u(h)-u(k)|^2 \Big)^{\frac{1}{2}}  \\&\geq \sum_{i \in \mathbb{Z}^d} u^2(i)- \frac{2\sqrt{d}}{\sqrt{d!}} \Big(\sum_{i \in \mathbb{Z}^d} u^2(i) \Big)^{\frac{1}{2}}  \Big(d!\sum_{(h,k) \in \mathcal{N}}|u(h)-u(k)|^2 \Big)^{\frac{1}{2}}\,,\nonumber
\end{align}
where we used that, due to Proposition~\ref{Kuhn_properties_PROP}(d),
\begin{align*}
\sum_{i \in \mathbb{Z}^d} \sum_{\pi \in \mathcal{P}_d} \underset{h,k\in \mathbb{Z}^d \cap T_\pi(i)}{\sum_{(h,k) \in \mathcal{N}}}|u(h)-u(k)|^2 &= \sum_{i \in \mathbb{Z}^d} \sum_{\pi \in \mathcal{P}_d}\sum_{(h,k) \in \mathcal{N}}\chi_{T_\pi(i)}(h)\cdot\chi_{T_\pi(i)}(k) |u(h)-u(k)|^2 \\& =\sum_{(h,k) \in \mathcal{N}}\sum_{i \in \mathbb{Z}^d} \sum_{\pi \in \mathcal{P}_d}\chi_{h-T_\pi}(i)\cdot\chi_{k-T_\pi}(i) |u(h)-u(k)|^2
\\& =\sum_{(h,k) \in \mathcal{N}}d!\,|u(h)-u(k)|^2\,.
\end{align*}
Then, we obtain the bound in {\rm (b)} from  \eqref{lemma:poincare}. 
\end{proof}
\subsection{Comparison between discrete and continuum eigenvalues of a set} \label{subsec:Eigenvalue_comparison_set}
In the following lemma, we give an upper bound for the continuum eigenvalue of $\zeta(X)$ in terms of the discrete eigenvalue $\lambda_N(X)$. 
\begin{lemma}\label{lemma_beta}
    Let $N\in\mathbb{N}$ and $X \in \mathcal{X}_N$.  Then, there exists $C_d>0$ such that
    \begin{equation*}
        \lambda\big(\zeta(X)\big) \leq N^{-\frac{2}{d}}\lambda_{N}\big(X\big) + C_d N^{-\frac{3}{d}} \,.
    \end{equation*}
\end{lemma}

\begin{proof} 
    For $X \in \mathcal{X}_N$ let $u\colon \mathbb{Z}^d \to \mathbb{R}$ be a minimizer of \eqref{defin_discrete_F}, that is 
    \begin{align*}
        N^{-1+\frac{2}{d}} \sum_{(i,j) \in \mathcal{N}} |u(i)-u(j)|^2 = \lambda_N(X)\,, \quad \frac{1}{N}\sum_{i \in \mathbb{Z}^d}u^2(i)=1 \,, \quad u(i) =0 \text{ for all } i \in \mathbb{Z}^d\setminus X\,.
    \end{align*}
     Note that the piecewise-affine interpolation $\hat{u}$  introduced in Section~\ref{sec:zetaandaffine} is a competitor for $\lambda(\zeta(X))$. Thus, the claim follows by estimating the Rayleigh quotient of $\hat{u}$  for $N$ large enough by means of Lemma~\ref{utilde_properties_LEMMA} as follows:
    \begin{align*}
        \frac{\int_{\mathbb{R}^d}|\nabla \hat{u}|^2\,\mathrm{d}x}{\int_{\mathbb{R}^d}|\hat{u}|^2\,\mathrm{d}x} &\leq \frac{ N^{\frac{d-2}{d}}E_N(u) }{N-2\sqrt{d}\sqrt{N}(N^{\frac{d-2}{d}}E_N(u))^{\frac{1}{2}}} = \frac{N^{\frac{d-2}{d}}\lambda_N(X)}{N} \frac{1}{1-2\sqrt{d}\left(\frac{N^{\frac{d-2}{d}}\lambda_N(X)}{N}\right)^{\frac{1}{2}}}\\
        &= N^{-\frac{2}{d}}\lambda_N(X)\frac{1}{1-2\sqrt{d} N^{-\frac{1}{d}}\lambda^{\frac{1}{2}}_N(X)} \,.
    \end{align*}
    Furthermore, we can use that $\frac{1}{1-t}\le 1+2t$ for $t\in \left]0,\frac{1}{2}\right[$ to estimate the last right hand side in the previous inequality for $\lambda_N(X)\leq (16d)^{-1}N^{\frac{2}{d}}$ and get
    \begin{align*}
        \frac{\int_{\mathbb{R}^d}|\nabla \hat{u}|^2\,\mathrm{d}x}{\int_{\mathbb{R}^d}|\hat{u}|^2\,\mathrm{d}x} \leq N^{-\frac{2}{d}}\lambda_N(X) \left(1+4 \sqrt{d} N^{-\frac{1}{d}}\lambda_N^{\frac{1}{2}}(X)\right) = N^{-\frac{2}{d}}\lambda_N(X) + 2\sqrt{d} N^{-\frac{3}{d}}\,.
    \end{align*}
The statement of the lemma follows by passing to the infimum over $H^1(\Rd)$ functions which are zero outside $\zeta(X)$ and satisfy $\norm{u}_{L^2(\Rd)}=N$.
\end{proof}
\subsection{Comparison between minimal eigenvalues in the discrete and continuum setting}\label{subsec:min_eigenvalue} 
In this section we compare the two minimal eigenvalues for fixed cardinality $N$ (resp.~measure in the continuum) defined in \eqref{defin_discrete_F} and \eqref{ineq:Faber_Krahn}.

\begin{lemma}\label{lemma_gamma} 
    Let $N \in \mathbb{N}$ and let $m_{\lambda,N}$ be as in \eqref{def:m_lambda_N}. Then, there exists $C_d>0$ such that
    \begin{equation*}
        m_{\lambda,N}\leq \lambda\big(B_1\big)+C_dN^{-\frac{1}{d}}\,.
    \end{equation*}
\end{lemma}

\begin{proof}  
    We first prove the claim for $N_k:= \card(\overline{B}_k \cap \mathbb{Z}^d)$, $k \in \mathbb{N}$. Note that there exists $C_d>0$ such that
    \begin{align}\label{ineq:comparison-k-N}
        |B_1|\left(k^d -C_d k^{d-1}\right) \leq N_k \leq |B_1|\left(k^d +C_d k^{d-1}\right) \,.
    \end{align}
    In the rest of the proof we omit the dependence of $N$ on $k$ to simplify notation. Let $u \in H^1_0(B_1)$ be such that
    \begin{align*}
        \int_{B_1} |\nabla u|^2\,\mathrm{d}x = \lambda(B_1)\,, \quad \int_{\mathbb{R}^d} u^2\,\mathrm{d}x=1\,.
    \end{align*}
    The following scaling property can be readily checked to hold true:
    \begin{align}\label{eq:scalinglambdacont}
        \lambda(B_1) = |B_1|^{-\frac{2}{d}}\lambda(B)\,,
    \end{align}
    where $B$ is the Euclidan Ball of unitary measure.
    By standard elliptic regularity estimates $u \in C^\infty(\overline{B}_1)$. Hence, for any $m \in \mathbb{N}$ and every multi-index $\alpha$ of order $m$ there exists $\overline{C}_{d,m}>0$ such that 
    \begin{align}\label{ineq:derivative-estimate}
        \|D^\alpha u\|_{L^\infty(\mathbb{R}^d)} \leq \overline{C}_{d,m}\,.
    \end{align}
    Let $u_k \colon \mathbb{Z}^d\to \mathbb{R}$ be defined as $u_k(i)= C_{d,k} u\left(\frac{i}{k}\right)$, $i \in \mathbb{Z}^d$, and where $C_{d,k}>0$ is to be chosen such that $\sum_{i \in \mathbb{Z}^d} u^2_k(i)=N$. Using \eqref{ineq:comparison-k-N} and the fact that $u \in H^1_0(B_1)$, we can compute $C_{d,k}$ as follows. We first observe that
    \begin{align*}
        \sum_{i \in \mathbb{Z}^d \cap \overline{B}_k} \left(u\left(\frac{i}{k}\right)\right)^2 &= \sum_{i \in \frac{1}{k}\mathbb{Z}^d \cap \overline{B}_1} u^2\left(i\right) \\&= \sum_{i \in \frac{1}{k}\mathbb{Z}^d \cap \overline{B}_1} \left( u^2\left(i\right)-\fint_{Q_{\frac{1}{k}}(i)}u^2(x)\,\mathrm{d}x\right) + k^d\int_{\mathbb{R}^d} u^2(x)\,\mathrm{d}x \\&= \sum_{i \in \frac{1}{k}\mathbb{Z}^d \cap \overline{B}_1} \left( u^2\left(i\right)-\fint_{Q_{\frac{1}{k}}(i)}u^2(x)\,\mathrm{d}x\right)  +\frac{N}{|B_1|} +O(k^{d-1})\,.
    \end{align*}
    By the regularity of $u$ in \eqref{ineq:derivative-estimate}, we can estimate
    \begin{equation*}
        \sum_{i \in \frac{1}{k}\mathbb{Z}^d \cap \overline{B}_1} \left| u^2\left(i\right)-\fint_{Q_{\frac{1}{k}}(i)}u^2(x)\,\mathrm{d}x\right| \leq C_d k^{d-1} \|u\|_{L^\infty} \cdot \norm{\nabla u}_{L^\infty} \leq C_d k^{d-1}\,.
    \end{equation*}
    This implies that $C_{d,k}$ satisfies
    \begin{equation}\label{ineq:Cdk}
        |B_1|^{1/2}\left(1-\frac{C_d}{k}\right)\leq C_{d,k}\leq |B_1|^{1/2}\left(1+\frac{C_d}{k}\right)\,.
    \end{equation}
    Since $u_k(i) =0$  for all $i \in \mathbb{Z}^d \setminus B_k$, $u_k$ is an admissible competitor for $\lambda_N(X)$ with $X \in \mathcal{X}_N$ given by $X= \overline{B}_k\cap \mathbb{Z}^d$. Therefore we have that
    \begin{align}\label{ineq:mlambdaN}
        m_{\lambda,N} \leq \lambda_N(X_{N}) \leq N^{-1+\frac{2}{d}} \sum_{(i,j) \in \mathcal{N}} |u_k(i)-u_k(j)|^2=E_{N}(u_k)\,.
    \end{align}
    As a consequence, in order to show the estimate for $m_{\lambda,N}$, it suffices to estimate $E_N(u_k)$. 
   To this end we note that,
%
    by the regularity of $u$ given in \eqref{ineq:derivative-estimate}, there exists $C_d >0$ such that
    \begin{align*}
     |u_k({i+e_n})-u_k(i)|^2 \leq C_{d,k}^2 k^{d-2} \int_{Q_{\frac{1}{k}}(\frac{i}{k})}  \left|\partial_n u\left(x\right)\right|^2 \,\mathrm{d}x + C_{d} k^{d-3} \left|Q_{\frac{1}{k}}\left(\frac{i}{k}\right)\right|
     \,.
    \end{align*}
    Summing over $i \in \mathbb{Z}^d$, $n \in \{1,\ldots,d\}$ and using \eqref{ineq:comparison-k-N},\eqref{eq:scalinglambdacont}, and \eqref{ineq:Cdk}, we obtain
    \begin{align*}
    E_N(u_k)=N^{-1+\frac{2}{d}}\sum_{(i,j)\in \mathcal{N}} |u_k({i})-u_k(j)|^2 &\leq    N^{-1+\frac{2}{d}} C_{d,k}^2 k^{d-2} \int_{B_1}  \left|\nabla u\left(x\right)\right|^2 \,\mathrm{d}x  + C_d N^{-1+\frac{2}{d}}C_{d,k}^2 k^{d-3}  \\&\leq  N^{-1+\frac{2}{d}} |B_1| k^{d-2} \int_{B_1}  \left|\nabla u\left(x\right)\right|^2 \,\mathrm{d}x  + C_d N^{-1+\frac{2}{d}} k^{d-3}  \\&\leq \lambda(B_1) + C_d N^{-\frac{1}{d}}\,.
    \end{align*}
    The latter estimate together with \eqref{ineq:mlambdaN} concludes the proof in the case that $N= \card(\overline{B}_k \cap \mathbb{Z}^d)$ for $k \in \mathbb{N}$. To treat the general case we note that for each $N \in \mathbb{N}$ there exists $k \in \mathbb{N}$ such that for $N_k = \card(\overline{B}_k \cap \mathbb{Z}^d)$ we have
    $
    0\leq  N-N_k\leq C_d k^{d-1}
    $. Finally, the estimate in the statement follows by using as a test the function $u_k$ constructed above and as $X_N=X_{N_k} \cup Z_N$ where $Z_N \subset \mathbb{Z}^d$ is chosen such that $\card X_N=N$.
\end{proof}

\subsection{Comparison of the asymmetries in the discrete and continuum setting} \label{subsec:Fraenkel}

In this section, we estimate the asymmetry of a set $X \in \mathcal{X}_N$ with respect to a (properly discretized) ball with that of its continuum embedding $\zeta(X)$.

\begin{lemma}\label{lemma_delta}    Let $N\in\mathbb{N}$ and $X \in \mathcal{X}_N$.  Then, there exists $C_d>0$ such that for all $z \in \mathbb{Z}^d$ there holds
    \begin{equation*}\label{discrete_asymmetry_inequality}
        \card \left( X\Delta( z+ B_{r_N}\cap\Z^d) \right) \leq \left|\zeta(X)\Delta (B_{r_N}+z)\right| +C_d N^{\frac{d-1}{d}} P_N(X) \,,
    \end{equation*}
    where $ r_N >0 $ is such that $|B_{r_N}|=N$.
\end{lemma}

\begin{proof} 
    Without loss of generality we assume $z =0$. In order to prove the statement, we claim that for $X,Y \subset \mathbb{Z}^d$ we have
    \begin{align}\label{ineq:symmetricdiffdiscretecont}
        \card\left( X \Delta Y \right) \leq  |\zeta(X) \Delta \zeta(Y)| + C_d \left(P(X) +P(Y)\right)\,.
    \end{align}
    Assuming the claim, we now show how to conclude. If $Y = B_{r_N} \cap \mathbb{Z}^d$, then $|\card Y-N|\leq C_d N^\frac{d-1}{d}$ which implies $P(Y)\leq C_d N^{\frac{d-1}{d}}$. Since $X \in \mathcal{X}_N$, $\card X=N$, hence, by the isoperimetric inequality on $\mathbb{Z}^d$, there exists $C_d>0$ such that $P(X)\geq C_d N^{\frac{d-1}{d}}$. Thus, we deduce that
    \begin{align*}
        P(Y) \leq C_d P(X)\,.
    \end{align*}
Arguing as in the proof of Lemma~\ref{lemma_measure_zeta}, for $Y = B_{r_N} \cap \mathbb{Z}^d$ it holds that
\begin{equation}\label{eq:zetaYBrN}
        |\zeta(Y) \Delta B_{r_N}| \leq  C_d P(Y)\leq C_d P(X)\,.
    \end{equation}
    Recalling Definition~\ref{def:Per}, the statement of the lemma eventually follows from \eqref{ineq:symmetricdiffdiscretecont}, \eqref{eq:zetaYBrN}, and the triangle inequality recalling the scaling $P(X) = N^{\frac{d-1}{d}}P_N(X)$. We now prove the claim \eqref{ineq:symmetricdiffdiscretecont}. We first observe that \eqref{ineq:symmetricdiffdiscretecont} follows from the estimate
    \begin{align}\label{ineq:XsetminusY}
        \card\left( X\setminus Y \right) \leq |\zeta(X) \setminus \zeta(Y)| +  C_dP(Y)
    \end{align}
    by exchanging the roles of $X$ and $Y$. The proof of \eqref{ineq:XsetminusY} is the consequence of the following two facts. On one hand, by the very definition of $\zeta$ in \eqref{def:zeta} it follows that
    \[
    \bigcup_{x\in X\setminus \overline{\zeta(Y)}} (x+Q_1) \subset \zeta(X)\setminus \zeta(Y)\,.
    \]
    Hence, we infer that
    \[
    \card \left(X\setminus \overline{\zeta(Y)}\right) = \sum_{x\in X\setminus \overline{\zeta(Y)}}|Q_1| = \left|\bigcup_{x\in X\setminus \overline{\zeta(Y)}} (x+Q_1)\right| \le |\zeta(X)\setminus \zeta(Y)|\,.
    \]
    On the other hand, since any $x\in \overline{\zeta(Y)}\setminus Y$ is a neighbour of a boundary point of $Y$, it follows that
    \begin{align*}
        \card\{x \in X \setminus Y \colon x\in\overline{\zeta(Y)}\} \leq C_d P(Y)\,.
    \end{align*}
    Combining these two inequalities, we get \eqref{ineq:XsetminusY} and conclude the proof.
\end{proof}

\EEE

\subsection{Proof of the main result}

\begin{proof}[Proof of Theorem~\ref{main_THM}]  Using Lemma~\ref{lemma_gamma} and \eqref{ineq:Xalmostmin}, we have that
\begin{align*}
    \lambda_N(X) \leq \lambda(B_1) + \alpha_N +C_d N^{-\frac{1}{d}}\,.
\end{align*}
From Lemma \ref{lemma_beta} and the above estimate we infer that
\begin{align*}
    \lambda(\zeta(X)) \leq N^{-\frac{2}{d}}\lambda(B_1) + N^{-\frac{2}{d}}\alpha_N  +C_d N^{-\frac{3}{d}} \,.
\end{align*}
Furthermore, setting $r_{|\zeta(X)|}>0$ such that $B_{r_{|\zeta(X)|}}=|\zeta(X)|$, by the scaling properties of $\lambda$, Lemma~\ref{lemma_measure_zeta} and the fact that $t\mapsto t^{\frac{2}{d}}$ is concave, we have
\begin{align*}
    N^{-\frac{2}{d}}\lambda(B_1) = \left( \frac{|\zeta(X)|}{N}\right)^{\frac{2}{d}} \lambda(B_{r_{|\zeta(X)|}})\leq \lambda(B_{r_{|\zeta(X)|}})\left(1+ \frac{2}{d}C_d N^{-\frac{1}{d}}P_N(X)\right) \,.
\end{align*}
 By Lemma~\ref{lemma_measure_zeta} we have $|\zeta(X)|\leq C_d N$ and therefore, again by the scaling properties of $\lambda$, we have
\begin{align*}
    \frac{\lambda(\zeta(X))-\lambda(B_{r_{|\zeta(X)|}}) }{\lambda(B_{r_{|\zeta(X)|}})}\leq C_d N^{-\frac{1}{d}} P_N(X) + C_d\alpha_N  +  C_d N^{-\frac{1}{d}} \leq  C_d N^{-\frac{1}{d}} P_N(X) + C_d\alpha_N\,.
\end{align*}
Using the subadditivity of the square-root we obtain
\begin{align*}
    \left(\frac{\lambda(\zeta(X))-\lambda(B_{r_{|\zeta(X)|}}) }{\lambda(B_{r_{|\zeta(X)|}})}\right)^{\frac{1}{2}} \leq  C_d\left( N^{-\frac{1}{2d}} P_N(X)^{\frac{1}{2}} + \alpha_N^{\frac{1}{2}}\right)\,.
\end{align*}
By \cite[Main Theorem]{BDPV} there exists $z \in \mathbb{R}^d$ such that 
\begin{equation*}
    \left|\zeta(X) \Delta( B_{r_{|\zeta(X)|}} +z)\right| \leq |\zeta(X)| C_d \left( N^{-\frac{1}{2d}} P_N(X)^{\frac{1}{2}} +\alpha_N^{\frac{1}{2}}\right)\,.
\end{equation*}  
Combining the latter inequality with the estimate in Lemma \ref{lemma_measure_zeta}, we infer that
\begin{equation}\label{ineq:FKpartial}
    \left|\zeta(X) \Delta( B_{r_{|\zeta(X)|}} +z)\right| \leq  C_dN \left( N^{-\frac{1}{2d}} P_N(X)^{\frac{1}{2}} +\alpha_N^{\frac{1}{2}}\right)\,.
\end{equation}
Let $\overline{z} \in \mathbb{Z}^d$ be such that $|z-\overline{z}|\leq \sqrt{d}$. Then, the triangle inequality gives that
\begin{align*}
    |(B_{r_N} +\overline{z})\Delta( B_{r_{|\zeta(X)|}} +z)| &\leq \left|(B_{r_N} +\overline{z})\Delta( B_{r_{|\zeta(X)|}} +\overline{z})\right| + \left|( B_{r_{|\zeta(X)|}} +\overline{z}) \Delta( B_{r_{|\zeta(X)|}} +z)\right| \,.
\end{align*}
Thanks to Lemma~\ref{lemma_measure_zeta}, the first term on the right hand side can be estimated as
\begin{equation*}
    \left|(B_{r_N} +\overline{z})\Delta( B_{r_{|\zeta(X)|}} +\overline{z})\right| \leq C_d N^{\frac{d-1}{d}}P_N(X)\,.
\end{equation*}
The second term on the right hand side can be estimated by evaluating the measure of the symmetric difference between a ball of radius $r_{|\zeta(X)|}+\sqrt{d}$ and a ball of radius $r_{|\zeta(X)|}-\sqrt{d}$ and then by using Lemma~\ref{lemma_measure_zeta}: 
\begin{align*}
    \left|( B_{r_{|\zeta(X)|}} +\overline{z}) \Delta( B_{r_{|\zeta(X)|}} +z)\right| &\leq |\zeta(X)|\left(\left(1+\frac{\sqrt{d}}{r_{|\zeta(X)|}}\right)^d-\left(1-\frac{\sqrt{d}}{r_{|\zeta(X)|}}\right)^d\right) \\
    &\leq C_d |\zeta(X)|\left(\frac{2d\sqrt{d}}{r_{|\zeta(X)|}}\right) \leq C_d \frac{|\zeta(X)|}{N^{\frac{1}{d}}}\\
    &\leq C_d N^{\frac{d-1}{d}}\left(1+C_dN^{-\frac{1}{d}}P_N(X)\right)\,.
\end{align*}
We can thus write that 
\begin{equation}\label{eq:balls_asymmetry}
    \left|(B_{r_N} +\overline{z})\Delta( B_{r_{|\zeta(X)|}} +z)\right| \leq  C_d N^{\frac{d-1}{d}}P_N(X)+ C_dN^{\frac{d-1}{d}}\,.
\end{equation}
Combining the estimate in Lemma \ref{lemma_delta} and \eqref{eq:balls_asymmetry}, we obtain
\begin{align*}
    \card\left(X\Delta( \mathbb{Z}^d \cap B_{r_N}+\overline{z})\right) &\leq |\zeta(X) \Delta(B_{r_N} +\overline{z})|+ C_d N^{\frac{d-1}{d}}P_N(X)\\
    &\leq \left|\zeta(X) \Delta( B_{r_{|\zeta(X)|}} +z)\right| +\left|(B_{r_N} +\overline{z})\Delta( B_{r_{|\zeta(X)|}} +z)\right| + C_d N^{\frac{d-1}{d}}P_N(X)\\
    &\leq \left|\zeta(X) \Delta( B_{r_{|\zeta(X)|}} +z)\right| + C_d N^{\frac{d-1}{d}}P_N(X)+ C_dN^{\frac{d-1}{d}}\,. 
\end{align*}
 Eventually, by \eqref{measure_zeta_estimate} and \eqref{ineq:FKpartial}, we infer that 
\begin{align*}
    \card\left(X\Delta( \mathbb{Z}^d \cap B_{r_N}+\overline{z}) \right) &\leq |\zeta(X) \Delta( B_{r_{|\zeta(X)|}} +z)| + C_d N^{\frac{d-1}{d}}P_N(X)+ C_dN^{\frac{d-1}{d}}\\
    &\leq  C_dN\left( N^{-\frac{1}{2d}} P_N(X)^{\frac{1}{2}} +\alpha_N^{\frac{1}{2}}  \right) +
    C_d N^{\frac{d-1}{d}}P_N(X)+ C_dN^{\frac{d-1}{d}} \\& \leq C_d N\left(
   N^{-\frac{1}{2d}} P_N(X)^{\frac{1}{2}}  + \alpha_N^{\frac{1}{2}} +N^{-\frac{1}{d}} \right)\,,
\end{align*}
where in the last inequality we used that $P_N(X)\leq C_d N^{\frac{1}{d}}$.
To conclude, we observe that $\sup_{N}\alpha_N <+\infty$ yields \eqref{ineq:fluctuation}, while \eqref{ineq:quasiminimi} follows from the assumptions on $X_N$ and the estimate~\eqref{ineq:fluctuation}. 
\end{proof}

%
%
%

\section*{Acknowledgments}
 The research of L.~Kreutz  was supported by the DFG through the Emmy Noether Programme (project number 509436910). The research of G.P.~Leonardi has been supported by Project PRIN 2022PJ9EFL ``Geometric Measure Theory: Structure of Singular Measures, Regularity Theory and Applications in the Calculus of Variations'' (financed by European Union - Next Generation EU, Mission 4, Component 2, CUP\_E53D23005860006). The research of G.~Morselli has been supported by INdAM -- GNAMPA Project 2023: ``Esistenza e propriet\`a fini di forme ottime'' (project CUP\_E53C22001930001). 

\EEE

\end{document}